\documentclass{article}
\usepackage[utf8]{inputenc}
\def\MR#1{\href{http://www.ams.org/mathscinet-getitem?mr=#1}{MR#1}}
\usepackage{amssymb}
\usepackage{hyperref}
\usepackage{amsmath}
\newtheorem{lem}{Lemma}
\newtheorem{ter}{Theorem}
\newtheorem{utv}{Proposition}
\newtheorem{sled}{Corollary}
\title{Flow polynomials as Feynman amplitudes and their $\alpha$-representation
}
\author{Eduard~Yu.~Lerner, Andrey~P.~Kuptsov, Sofya~A.~Mukhamedjanova}
\date{}

\begin{document}
\maketitle
\begin{abstract}       
  Let $G$ be a connected graph; denote by $\tau(G)$ the set of its spanning
  trees. Let $\mathbb F_q$ be a finite field, $s(\alpha,G)=\sum_{T\in\tau(G)}
  \prod_{e \in E(T)} \alpha_e$, where ${\alpha_e\in \mathbb F_q}$. Kontsevich
  conjectured in 1997 that the number of nonzero values of $s(\alpha, G)$ is
  a polynomial in $q$ for all graphs. This conjecture was disproved by
  Brosnan and Belkale. In this paper, using the standard technique
  of the Fourier transformation of Feynman amplitudes, we express the flow
  polynomial $F_G(q)$ in terms of the ``correct'' Kontsevich formula.
  Our formula represents~$F_G(q)$ as a linear combination of Legendre
  symbols of $s(\alpha, H)$ with coefficients $\pm 1/q^{(|V(H)|-1)/2}$, where~$H$
  is a contracted graph of~$G$ depending on~$\alpha\in \left(\mathbb F^*_q\right)^{E(G)}$,
  and $|V(H)|$ is odd. The case $q=5$ corresponds to the least number with which
  all coefficients in the linear combination are positive. This allows us to hope
  that the obtained result can be applied to prove the Tutte 5-flow conjecture.
\end{abstract}


\section{Introduction. The statement of the main result}

Let $q$ be a positive integer and let $A_q$ be an arbitrary abelian group consisting of $q$ elements; we usually use the additive group of the field ${\mathbb F}_q$ for $A_q$; in this case $q=p^d$, 
where $p$ is prime and $d\in\mathbb N$.
Let $G(V,E)$ be a connected multigraph without loops; let $V(G)$ denote the set of its vertices and $E(G)$ do the set of its edges. When considering the initial graph, we sometimes omit the symbol~$G$ in denotations. In certain cases we need to indicate the orientation of graph edges, so we denote the origin of an edge $e\in E(G)$ by $i(e)$ and do its terminus by $f(e)$.

Recall that the chromatic polynomial $P_G(q)$ counts the number of proper vertex colorings with~$q$ colors. Define the norm of an element~$y$ of the group~$A_q$ as
\[
||y|| = \left\{
\begin{array}{ll} 0, & y = 0,  \\
1, & y \ne 0,
\end{array}
\right.
\]
where the symbol $0$ in the right-hand side of the equality $y=0$ is the neutral element of the group. Let us associate each vertex~$v$ with its color~$x_v$. Evidently, a coloring is proper if and only if
\[
\prod_{e\in E(G)} || x_{i(e)}  - x_{f(e)}||=1,
\]
and otherwise the product equals 0.

Using this fact, we get
\[
  P_G (q) = \sum_{\begin{subarray}{e}
  x_v  \in A_q  \\ \forall v \in V(G)\end{subarray} }
  \prod_{e\in E(G)} || x_{i(e)}  - x_{f(e)}||.
\]
This representation allows us to treat chromatic polynomials as vacuum Feynman amplitudes in the coordinate space (see Section~2 for details).

Let us define the delta function on an abelian group ($x\in A_q$) by the formula
$$
\delta (x) = \left\{ {\begin{array}{ll} 0, & {x \ne 0},  \\
1, & {x = 0},  \end{array} } \right.
$$
i.e., in our case $\delta(x)=1-||x||$. Let $(\varepsilon_{v
e})_{v\in V,e\in E}$ be the incidence matrix of an arbitrary orientation of the graph $G$; it obeys the formula
$$
\varepsilon_{v e}  = \left\{
  \begin{array}{ll}
   -1, & \mbox{if $i(e)=v$},  \\
   1, & \mbox{if $f(e)=v$},  \\
   0, & \mbox{if $e$ is nonincident to $v$}.
  \end{array} \right.
$$

Let us associate each edge~$e$ of the graph~$G$ with an element $k_e$ of the group~$A_q$ so that $\sum_{e\in E}\varepsilon_{v e} k_e=0$. This association is called a flow. The number of everywhere nonzero flows is a polynomial in~$q$; it is called the flow polynomial~\cite{distel,Tutte}. Therefore, the flow polynomial obeys the formula
\begin{equation} \label{eq:1}
  F_G(q) = \sum_{
  \begin{subarray}{e} k_e  \in A^*_q  \\ \forall e\in E(G) \end{subarray}}
  {\prod_{v \in V(G)} {\delta \left( {\sum_{e\in E(G)}
  {\varepsilon_{v e} } k_e } \right)} },
\end{equation}
(here $A_q^*$ is the collection of nonzero elements of the group~$A_q$). Another variant of formula~(\ref{eq:1}) takes the form
\begin{equation} \label{FGq}
  F_G(q) = \sum_{
  \begin{subarray}{e} k_e  \in A_q  \\ \forall e\in E(G) \end{subarray}}
 \prod_{e \in E(G)} ||k_e|| \prod_{v \in V(G)} \delta \left( {\sum_{e\in E(G)}
  \varepsilon_{v e}  k_e } \right).
\end{equation}
Formula~(\ref{FGq}) allows us to treat flow polynomials as vacuum Feynman amplitudes in the impulse space (see Section~2). Using the technique of Feynman amplitudes, we can get a new representation for the flow polynomial. Let us now state the main result of this paper.

The main theorem considers the case when $q = p^d$, $p$ is an odd prime, and $d\in\mathbb N$. Denote by $\eta$ the multiplicative quadratic character of the field~$\mathbb{F}_q$: $\eta(0)=0$, in other cases $\eta(x)=1$ or $\eta(x)=-1$ in dependence of whether $x$ is a square in the field $\mathbb F_q$ or not. For $d=1$ the function $\eta$ coincides with the Legendre symbol of the residue field modulo prime~$p$. Denote by $g(q)$ the quadratic Gaussian sum of the field~$\mathbb{F}_q$. It obeys the formula
\begin{equation}
\label{gq}
g(q) = \left\{ \begin{array}{ll}
   (-1)^{d-1}\sqrt{q},& \textrm{if}~p \bmod 4= 1,  \\
   (-1)^{d-1} i^d \sqrt{q},&  \textrm{if}~p \bmod 4=3,  \end{array} \right.
\end{equation}
here $i$ is the imaginary unit.

Let all edges of the graph~$G$ be associated with nonzero weights $\alpha_e$: $\alpha_e\in {\mathbb F}_q^*$ (${\mathbb F}_q^*$ is the collection of nonzero elements of the field). Denote the graph with the contracted vertex set $W$ ($W\subseteq V$) by $G/W$; note that it contains no internal edges of the subgraph~$H(W)$. Weights $\alpha_e$ of edges of the graph~$G/W$ are equal to weights of corresponding edges of the graph~$G$.

Let $\tau(G)$ be the set of spanning trees of the graph~$G$. Consider sums
\begin{equation} \label{eq:2}
s(\alpha,G)=\sum_{T \in \tau(G)}
  {\prod_{e\in E(T)} {\alpha_e } },
\end{equation}
which, evidently, are elements of the field ${\mathbb F}_q$; if the graph~$G$ consists of one vertex, then by definition we put $s(\alpha,G)=1$.

Denote by $W^*=W^*(\alpha,G)$ any nonempty minimum cardinality subset $W$ of the vertex set $V(G)$, for which the sum $s(\alpha,G/W)$ differs from zero. Evidently, if $|W|=1$ then $G/W\equiv G$. Note also that the number of edges in each tree $T\in\tau(G/W^*)$ equals $|V|-|W^*|$. Denote this difference by~$r^*(G,\alpha)$.
\begin{ter}[The main theorem]
\label{ter:1}
Let $G$ be an arbitrary connected multigraph and let $q = p^d$ with odd prime $p$ and $d\in\mathbb N$. Then
\begin{equation} \label{eq:3}
F_G (q) = \sum_{\alpha \in \left(\mathbb F^*_q\right)^{E(G)}} {\eta
\left(s(\alpha,G/W^*)\right) \left[ \frac{g(q)}q
\right]^{r^*(G,\alpha)} }.
\end{equation}
\end{ter}

Note that $s(\alpha,G)$ is an algebraic complement of any element of the weighted Laplacian matrix~$L$ of the graph $G$. As we prove in Section~4, $s(\alpha,G/W^*)$ coincides with the largest dimension nondegenerate principal minor of the matrix~$L$. Therefore, one can interpret our theorem as a representation of the flow polynomial as a linear combination of Legendre symbols of minors of the Laplacian matrix.

Note that in what follows we assume that the multigraph has no loops, though this is not explicitly stated in Theorem~\ref{ter:1}. The fact is (as one can easily see) that by removing a loop from the graph we make the left- and right-hand sides of equality~(\ref{eq:3}) exactly $q-1$ times less. So the proof of the theorem for an arbitrary multigraph is reduced to the proof for a multigraph without loops.

Consider the graph~$K_3$ as a simplest illustration of Theorem~\ref{ter:1}. Its flow polynomial obeys the formula $F_{K_3}(q)=q-1$. One can prove (see Propositions~\ref{utv:1}--\ref{utv:2} in Section~7) that in sum~(\ref{eq:3}) the terms that correspond to the case $r^*(G,\alpha)=1$ cancel each other out, and this sum contains no terms that correspond to the case $r^*(G,\alpha)=0$. Therefore, the assertion of Theorem~\ref{ter:1} for~$K_3$ is equivalent to the equality
\begin{equation}
\label{K3}
q-1 = \sum_{\alpha_1,\alpha_2,\alpha_3 \in \mathbb F^*_q} \eta(\alpha_1\alpha_2+\alpha_1\alpha_3+\alpha_2\alpha_3)
 \left[ \frac{g(q)}q\right]^2.
\end{equation}

A natural generalization of the notion of a flow polynomial for the case of an arbitrary matroid is the notion of the characteristic function of the dual matroid~\cite{aigner}. We discuss the generalization of Theorem~\ref{ter:1} 
for arbitrary matroids representable over the field~$\mathbb F_q$
in a separate publication.

The paper has the following structure. In Section~2 we give a brief information on Feynman amplitudes and motivate our interest to flow polynomials. This section is not necessary for a formal understanding of the proof of Theorem~\ref{ter:1}, but it is useful for the comprehension of its sources and of prospects for the use of the mentioned technique. In Section~3 we recall some properties of the Fourier transformation over a finite field and prove the key lemma which represents $F_G (q)$ as a double sum, namely, the sum that is first taken over $\alpha$ and then over~$k$. In Section~4 we prove that the matrix of the quadratic form in the exponent of the function summed with respect to~$k$ is the Laplacian matrix and recall its combinatorial properties related to the evaluation of its minors. In Section~5 we learn to calculate multidimensional Gaussian sums over a finite field. Finally, in Section~6 by using the obtained results we prove Theorem~\ref{ter:1}. Note that the proof is one paragraph long. In Section~7 we simplify formula~(\ref{eq:3}), using the fact that, in particular, its terms that correspond to odd values of $r^*(G,\alpha)$ cancel each other out. We also discuss possible applications of this formula for proving the Tutte 5-flow conjecture and give results of calculations for some graphs with~$q=5$.

\section{Some properties of Feynman amplitudes}

The simplest case when we are faced with vacuum Feynman amplitudes (FA) consists in the calculation of mean values
\begin{equation}
\langle \exp(-\sum_{x=1}^n \epsilon\, \varphi^4(x)) \rangle_{\mu_0},
\label{phi4}
\end{equation}
where $\mu_0$ is the Gaussian measure with the binary correlation function (the so-called propagator) $\rho(x,y)$, in formula~(\ref{phi4}) $x,y\in\{1,\ldots,n\}$ and $\varphi(x)\in \mathbb R$. As is known (see, for example, \cite[section~2.2]{malyshev}), one can use the so-called pairing technique for calculating the mean value with respect to the Gaussian measure, for example, $\langle \varphi^2 (x) \varphi^2 (y) \rangle_{\mu_0}=\rho(x,x)\rho(y,y)+2 \rho(x,y)^2$ (just for this reason $2n$th moments of the standard Gaussian distribution equal $(2n-1)!!$).

Expanding the exponent in~(\ref{phi4}) in a series, in the $n$th order with respect to $\epsilon$ we get homogeneous graphs with~$n$ vertices of the degree~4 labeled with variables like~$x,y$ (over which we calculate the sum), whose edges contribute the value $\rho(x,y)$ to the product of propagators (which is to be summed up).

Actually, Feynman used this technique in the case of a functional measure, when the index~$x$ itself took on values in a continual set like~${\mathbb R}^4$ rather than in a finite one, and instead of the summation with respect to $x$ there was the integration. In the so-called scalar models in the quantum field theory~\cite{smirnov} the propagator (in the coordinate representation) is set to $||x-y||^\lambda$ (here $x,y\in {\mathbb R}^4$). We can get the definition of $P_G(q)$ proposed by us above by exactly transferring all definitions given for variables that take on values in~${\mathbb R}^4$ to the case of a finite group (a finite field).

Note that we do not consider results of the evaluation of the vacuum variant of real FA with the propagator $||x-y||^\lambda$ because the integral diverges with all~$\lambda$ for any graph~$G$. It seems to be possible to avoid this obstacle by performing the integration in all variables, except one. Namely, since $\rho(x,y)$ depends only on the difference $x-y$, the result of the integration in $|V|-1$ variables is independent of the value of the rest one. Just such an integral is called a vacuum FA in the case of other transitive invariant propagators in the coordinate space, which explains the used terminology. However, actually, in the quantum field theory with the propagator $||x-y||^\lambda$ considered here the vacuum integral is always set to infinity.

Nonvacuum FA with such propagator are less trivial. Several variables (more than one) are fixed there, and the integration is performed over all the rest ones. An analog of a nonvacuum FA for a finite field gives the number of proper colorings of the graph~$G$, provided that some of its vertices have got colors $y_v$ already. Note also that the graph precoloring extension problem is well known~\cite{precoloring}, its complexity for various types of graphs is studied rather thoroughly. However, properties of the polynomial that define the number of proper colorings in this case are studied less.

We need these polynomials and their flow analogs for deducing explicit formulas for FA of a~$p$-adic argument; see section~3 in the paper~\cite{lerner} for the corresponding explicit formulas that take into account the specificity of the transfer of the known properties of chromatic and flow polynomials to the mentioned case. Note that here, on the contrary, we consider the application of the FA technique in combinatorics.

FA in the impulse space have been used from the very beginning of the development of this technique. If in the coordinate representation we consider a propagator in the form $\rho(x,y)=f(x-y)$, then in the impulse representation each edge is associated with the function $\widehat f(k)$, where $\widehat f$ is the Fourier transformation of~$f$. Since in ${\mathbb R}^4$ it holds $\widehat{||\cdot||^\lambda}=c(\lambda)
||\cdot||^{-4-\lambda}$ (one can easily verify this property; note that the Fourier transformation in the real-valued case is understood in the sense of generalized functions, though in what follows such details are inessential), each edge in the diagram of the impulse representation in the real-valued case also corresponds to some degree of the norm. Formula~(\ref{FGq}) is a formal calque of definitions of vacuum FA accepted in the so-called real scalar theory (in view of remarks on the convergence analogous to those given above for the coordinate space).

In a nonvacuum case, in the impulse representation some variables (some real-valued analogs of variables~$k_e$ introduced by us) are fixed and the result of the integration with respect to all the rest variables depends on them. 
The fact that nonvacuum FA in the impulse and coordinate representations are connected with each other by the Fourier transformation explains the terms ``impulse/coordinate representation''. In Lemma~\ref{lem:1} given below this connection is considered in the simplest case of finite fields.

Note that in the theory of FA it has long been discovered that when formally associating one and the same propagators with edges of the graph of an FA in coordinate and impulse spaces, the vacuum amplitude in the coordinate space for the planar graph~$G$ coincides with the corresponding amplitude in the impulse space for the dual graph~$\widetilde G$~\cite{bleher}. But if propagators are connected with each other by the Fourier transformation, then vacuum amplitudes (if they are defined) in coordinate and impulse representations coincide.

In combinatorics the coincidence of chromatic and flow polynomials of dual graphs has also long been known, namely, from the very inception of the concept of a flow polynomial. However, in a finite field the Fourier transformation of a norm is not a norm. For this reason the connection between flow and chromatic polynomials of one and the same graph, which is a result of the Fourier transformation of a norm, is more complex. We discuss this connection in a separate paper.

Let us return to real FA. If all propagators have the same degrees, then even in the nonvacuum case there arises a difficulty with the convergence of all finite-dimensional integrals that define amplitudes. This difficulty can be eliminated by labeling each edge with ``its own'' complex-valued degree of the propagator and the subsequent analytic continuation of the integration result. This approach implies the estimation of the convergence domain in the space $\lambda_e$, which is also problematic. Moreover, in the quantum field theory, it is important to be able to find poles of the analytic continuation. This can be done with the help of the so-called $\alpha$-representation.

In the theory with a homogeneous propagator, the technique of the $\alpha$-representation consists in the following steps. One represents each function in the form $||k_e||^{-4-\lambda}$ associated with the edge $e$ in the impulse space as $||k_e^2||^{(-4-\lambda)/2}$ and replace it with the value of the Fourier transformation result of the norm of $\alpha_e$ raised to the corresponding degree at the point $k_e^2$. When calculating an FA in the impulse space, changing the integration order (we first integrate with respect to $k$ and then do with respect to $\alpha$), we get a multidimensional Gaussian integral with respect to variables~$k_e$, $e\in E$ with a special matrix depending on $\alpha$ (the weighted Laplacian matrix of the graph~$G$).
The Gaussian integral equals the square root of the matrix determinant; the latter obeys formula~(\ref{eq:2}) (for real $\alpha$). As a result we get a notation for the FA in the impulse space as the integral with respect to $\alpha$ of a function whose behavior can be studied easily (just this is called the $\alpha$-representation).

Formula~(\ref{eq:3}) gives an analogous representation (in the vacuum variant) for $F_G(q)$ in the case of a finite field $\mathbb F_q$. Instead of the integral we get the sum over all nonzero values of $\alpha_e$ in ${\mathbb F}_q$, and the result is expressed via the Legendre symbol of minors of the same Laplacian matrix as in the real-valued case.

The $\alpha$-representation of FA was also mentioned in earlier papers in combinatorics. In December 1997, when giving a talk at the Gelfand seminar in Rutgers University, Maxim Kontsevich proposed a conjecture that for any connected multigraph~$G$ the number $N(G,q)$ of nonzero values of $s(\alpha,G)$ for $\alpha  \in \left(\mathbb{F}_q\right)^{E(G)}$ is a polynomial in~$q$. Though the conjecture was never published, it has aroused the interest of experts in combinatorics (see~\cite{stanleyArticle,chung}). Against expectations, sometime later this conjecture was refuted in a nonconstructive way \cite{belk}. Constructive examples of graphs, for which the conjecture is not valid, start to be found comparatively recently~\cite{Schnetz}.

Note that if formula~(\ref{eq:3}) contains only those terms that correspond to the maximal value of the rank $r^*(G,\alpha)$ of the Laplacian matrix of the graph~$G$ (i.e., $r^*(G,\alpha)=|V|-1$), then (accurate to the coefficient) it allows the representation
$$
\sum\limits_{ \alpha \in \left(
\mathbb{F}_q^*\right)^{E(G)}} \eta \left(s(\alpha,G)\right).
$$
We need to impose no additional constraint on the rank in the sum, because terms such that $s(\alpha,G)=0$ contribute nothing to the sum. Note also that the value $N(G,q)$ mentioned in the Kontsevich conjecture is, evidently, representable as
$$
N(G,q)=\sum\limits_{ \alpha \in \left( \mathbb{F}_q\right)^{E(G)}}
\eta^2 \left(s(\alpha,G)\right).
$$

The main result of this paper (the $\alpha$-representation) is a ``proper'' notation of the Kontsevich conjecture.
The linear combination obtained with the help of the technique of FA is really a polynomial in~$q$, more precisely, it is 
the flow polynomial of $G$. This polynomial is a source of many unsolved questions ``dual'' to the map coloring problem.

\section{The Fourier transformation and flow polynomials}

The Fourier transformation in the finite group $A_q$ is defined with the help of the notion of an additive character $\chi(x)$, $x\in A_q$. Recall that an additive character~\cite[chapter 5]{lidlNider} is a complex-valued function $\chi(x)$, $x\in A_q$, such that $\chi(x+y)= \chi(x)\chi(y)$ for any $x,y\in A_q$. Evidently, for the neutral element of the group it holds $\chi(0)=1$, therefore the value of the complex module of the character identically equals one. Evidently, for the neutral group element it holds $\chi(0)=1$, therefore $|\chi(x)|=1$ for any $x\in A_q$.

The character that identically equals one is said to be trivial. One can easily prove (\cite[theorem~5.4]{lidlNider}) that for any nontrivial character it holds
\begin{equation}
\label{chizerro}
\sum_{x\in A_q} \chi(x)=0.
\end{equation}

In what follows, for $A_q$ we choose only an additive group of the finite field $\mathbb F_q$, $q=p^d$. It is well known that (\cite[theorem~5.7]{lidlNider}) any character of this group $\mathbb{F}_q$ takes the form $\chi_k(x)=\chi_1(kx)$, where $k\in$ $\mathbb{F}_q$ and $\chi_1(x)=\exp{(2\pi i\, \operatorname{Tr}(x)/p)}$, while $\operatorname{Tr}(x)=x+x^p+x^{p^2}\ldots+x^{p^{d-1}}$. One can easily prove that formula~(\ref{chizerro}) in this case allows the form
\begin{equation}
\label{eq:5} \sum_{k\in \mathbb F_q}{\chi_1(k t)}=q\,\delta(t).
\end{equation}

For any function $f(x)$ whose argument $x$ takes on values in~$\mathbb F_q$ we define the \emph{Fourier transformation} $\widehat{f}(k)$, $k\in \mathbb F_q$, as
$$
\widehat{f}(k)=\sum_{x\in
\mathbb F_q}{f(x)\chi_1( k x)}.
$$
Formula~(\ref{eq:5}) easily implies the equality
$$
f(x)=\frac1q\sum_{k\in \mathbb F_q}\widehat f(k)\chi_1(- k x)
$$
(the inverse Fourier transformation formula). Note that formula~(\ref{eq:5}) means that (accurate to the multiplier) the delta-function and the unit are connected with each other by the Fourier transformation.

\begin{lem}
\label{lem:1} Let $G$ be a multigraph without loops. Then the product of characters has the following property:
\begin{equation}
\label{eq:6}
\sum_{x \in \left(\mathbb F_q\right)^{V(G)}} \prod_{e\in E(G)} \chi_1((x_{i(e)}
- x_{f(e)})k_e)   = \prod_{v \in V(G)} q\,\delta \left(
\sum_{e\in E(G)} {\varepsilon_{v e} } k_e \right).
\end{equation}
\end{lem}
\textbf{Proof:}
Evidently,
$$
\prod_{e\in E(G)} {\chi_1( (x_{i(e)}  - x_{f(e)} )k_e )}= \chi_1
(\sum_{e\in E(G)} (x_{i(e)}  - x_{f(e)} ) k_e ).
$$
Note that the argument of the additive character is representable as
$$
\sum_{e\in E(G)} {(x_{i(e)}  - x_{f(e)} )k_e }=
\sum_{v \in V(G)} {x_v \sum_{e\in E(G)}{\varepsilon
_{v e} } k_e }.
$$
Certain insignificant transformations give
$$
\sum_{\begin{subarray}{e}
x_v  \in \mathbb F_q\\
\forall v \in V(G) \end{subarray}}
\prod_{e\in E(G)} {\chi_1((x_{i(e)}-x_{f(e)} )k_e )}=
\sum_{\begin{subarray}{e} x_v  \in \mathbb F_q
\\ \forall v \in V(G) \end{subarray}}
\prod_{v \in V(G)} {\chi_1(x_v \sum_{e\in E(G)} {\varepsilon_{v e} }
k_e } ) =
$$
$$
= \prod_{v \in V(G)}\left( {\sum_{x_v\in \mathbb F_q} {\chi_1(x_v
\sum_{e\in E(G)} {\varepsilon_{v e} k_e }) } }\right).
$$
Applying formula~(\ref{eq:5}), we get
$$
\prod_{v \in V(G)} {\sum_{x_v\in \mathbb F_q} {\chi_1(x_v}
\sum_{e\in E(G)} {\varepsilon_{v e}  k_e } ) }  =
\prod_{v \in V(G)} {q\,\delta \left( {\sum_{e\in E(G)}
{\varepsilon_{v e} } k_e } \right)}. 
$$
$\square$

We can interpret Lemma~\ref{lem:1} as follows. Consider an FA in the coordinate space with the propagator $\delta(x-y)$ with ``external variables'' $z_e$, i.e.,
$$\sum_{x\in (\mathbb F_q)^{V(G)}}
\prod_{e \in E(G)} \delta(x_{i(e)}-x_{f(e)}-z_e).$$
Then the Fourier transformation (with respect to variables $z_e$, $e\in E$) coincides (accurate to a constant coefficient) with the FA in the impulse representation with the unit propagator.
\begin{lem}[The Key lemma]
\label{lem:2} Let $G$ be a connected multigraph without loops. Then
\begin{equation}
\label{eq:7}
F_G (q) = q^{ - |V(G)|}
\sum_{\alpha  \in (\mathbb F_q^*)^{E(G)}}
\sum_{ x_v  \in \mathbb (F_q)^{V(G)}}
\chi\left( \sum_{e\in E(G)} (x_{i(e)}  - x_{f(e)} )^2
\alpha_e  \right),
\end{equation}
\end{lem}
\textbf{Proof:}
Let us apply Lemma~\ref{lem:1} for $k_e=\alpha_e$ and calculate the sum over all nonzero $\alpha_e$. We get
\begin{equation}
\label{eq:8}
\frac
{\sum\limits_{\begin{subarray}{e} \alpha_e  \in
\mathbb F_q^*
\\  \forall e\in E(G)\end{subarray}}
\sum\limits_{\begin{subarray}{e} x_v  \in \mathbb F_q \\
\forall v \in V(G)\end{subarray}}
\prod\limits_{e\in E(G)}
\chi_1((x_{i(e)} - x_{f(e)} )\alpha_e)
}{q^{|V(G)|}} =
\sum_{
\begin{subarray}{e} \alpha_e  \in \mathbb F_q^*
\\  \forall e\in E(G)\end{subarray}}\prod_{v \in V(G)}
\delta \left( {\sum_{e\in E(G)} {\varepsilon_{v e} } \alpha_e
} \right).
\end{equation}
By the definition of a flow polynomial~(\ref{eq:1}) the right-hand side of the latter equality coincides with~$F_G(q)$.

Let us now change the summation order in the left-hand side of formula~(\ref{eq:8}). We get
$$
\sum_{\begin{subarray}{e} x_v  \in \mathbb F_q \\
\forall v \in V(G)\end{subarray}} {\prod_{e\in E(G)}
\sum_{\alpha_e  \in \mathbb F_q^*} {\chi_1((x_{i(e)} - x_{f(e)}
)\alpha_e) } }.
$$

One can easily see that
\begin{equation}
\label{eq:9}
\sum_{\alpha_e  \in \mathbb F_q^*}
\chi_1((x_{i(e)} - x_{f(e)} )\alpha_e)
 = \left\{
\begin{array}{ll}
q-1, & \text{if } x_{i(e)} = x_{f(e)},\\
-1, & \text{if }x_{i(e)} \ne x_{f(e)}.
\end{array}
\right.
\end{equation}
Really, let $x_{i(e)} - x_{f(e)} = y_e$. In view of~(\ref{eq:5}) with $y_e \ne 0 $ we get $\sum_{\alpha_e \in \mathbb F_q}
\chi_1(y_e \alpha_e)= 0$, otherwise $\sum_{\alpha_e \in \mathbb F_q}
\chi_1(y_e \alpha_e)= q$. In~(\ref{eq:9}) we calculate the sum over all $\alpha_e$, except $\alpha_e=0$, which corresponds to the term that equals~1.

The right-hand side of formula~(\ref{eq:9}) is a function of $||y_e|| =h(||y_e||)$, namely,
$$
h(z) = \left\{
{\begin{array}{ll} q-1,& \text{if }z=0,\\
-1, &\text{otherwise.} \end{array} } \right.$$

This is the key moment in our proof. Using this fact, we replace $x_{i(e)} - x_{f(e)}$ in the left-hand side of~(\ref{eq:8}) with $(x_{i(e)} - x_{f(e)})^2 $, while the value $||x_{i(e)}-x_{f(e)}||$ remains the same. We get the assertion of the lemma.
$\square$

\section {The matrix tree theorem}

In Lemma~\ref{lem:2} we have represented a flow polynomial as the sum of characters of a quadratic form with respect to variables $x$, namely, $$\sum_{e\in E(G)} (x_{i(e)} - x_{f(e)})^2 \alpha_e .$$ The matrix of this quadratic form is the weighted Laplacian matrix of the graph~$G$.
\begin{lem}
\label{lem:3} Let $G$ be a multigraph without loops. Then the following correlation takes place:
$$
  \sum_{e\in E(G)} ( x_{i(e)}  - x_{f(e)} )^2 \alpha _e  = x_V^t L x_V,
$$
where $x_V$ is the vector column of all variables associated with vertices; the superscript $t$ is the transposition sign; $L$ is the so-called weighted Laplacian matrix of the graph~$G$, i.e.,
\begin{equation}
\label{laplas} \ell_{kj}  = \left\{ \begin{array}{ll}
    - \sum_{e:\,\{i(e),f(e)\}=\{ k,j\} } \alpha _e,& k\ne j,  \\
   \sum_{e:\, k\in\{i(e),f(e)\} } \alpha _e ,&  k = j.  \end{array} \right.
\end{equation}
\end{lem}

Formula~~(\ref{laplas}) means that each nondiagonal element of the Laplacian matrix equals the sum of weights of \textit{all} edges connecting the corresponding vertices multiplied by $(-1)$; each diagonal element equals the sum of weights of all edges incident to the corresponding vertex. Therefore, the Laplacian matrix is symmetric and degenerate, the sum of elements in each row of this matrix equals zero.

Since Lemma~\ref{lem:3} is well known, we do not give its proof here. The case of a simple graph with $\alpha_e\equiv 1$ is studied in~\cite[lemma~4.3]{MatricesAndGraphs}. In what follows, in order to indicate the dependence of the matrix $L$ on $\alpha_e$, $e \in E(G)$, we denote this matrix by $L(G,\alpha)$; its determinant equals zero. Below we also need principal minors of the matrix $L(G,\alpha )$ of lesser orders. Let us discuss their combinatorial sense.

Let us first consider minors of the order~${|V|-1}$. This result is classical; it goes back to works by R.Kirchhoff, J.J.Sylvester, and A.Cayley published in the middle of the 19th century (see details in \cite[section~5.6, remarks to chapter~5]{stanleyBook},~\cite{stanleyArticle}). Recall that the symbol $s(\alpha,G)$ denotes sum~(\ref{eq:2}).

\begin{ter}[the matrix tree theorem,~\cite{stanleyBook, MatricesAndGraphs}]
\label{ter:2}  Let $G$ be a connected multigraph without loops, let $L'(G,\alpha )$ be obtained from $L(G,\alpha )$ by deleting the $i$th row and the $i$th column, $i\in V(G)$. Then with any $i$ it holds
$$
  \det L'(G,\alpha ) = s(\alpha,G),
$$
where the sum $s(\alpha,G)$ obeys formula~(\ref{eq:2}).
\end{ter}

Let us now consider minors of lesser orders. We are going to reduce this case to the previous theorem.
\begin{ter}
\label{ter:3} Let~$G$ be a connected multigraph without loops. The principal minor of the matrix $L(G,\alpha)$ that is formed by rows and columns with numbers \{$i_1,\dots i_k$\} coincides with $s(\alpha,G/W)$, where $W=V\setminus\{i_1,\dots i_k\}$. \quad
\end{ter}
\textbf{Proof:}
Really, all vertices of the graph $G/W$, except one ``contracted'' vertex~$v'$, correspond to vertices $i_1,\dots i_k$ of the initial graph. The submatrix of the matrix $L(G,\alpha)$ that is formed by rows $i_1,\dots i_k$ coincides with the submatrix obtained by deleting from the matrix $L(G/W,\alpha)$ the row and the column that correspond to the vertex $v'$. Applying the previous theorem for the graph $G/W$, we get the assertion of Theorem~\ref{ter:3}.
$\square$

Note that the idea of using trees of the graph $G/W$ for calculating minors
of the Laplacian matrix of the graph~$G$ is also not new,
it goes back to works of Alexander~Kelmans~\cite{kelmans}.
The following proposition is also an analog of Theorem~\ref{ter:3}
(apparently, it goes back to works of Miroslav~Fiedler).

\begin{ter}[see \cite{MatricesAndGraphs}, theorem~4.7]
\label{minor} Let~$G$ be a connected multigraph without loops. Then the principal minor of the matrix $L(G,\alpha)$ that is formed by deleting rows and columns with numbers ${i_1,\ldots i_k}$ equals
$$
\sum_{\{T_1,\ldots,T_k\} }\ 
  \prod_{e\in \cup_j E(T_j)} {\alpha _e } ;
$$
here the sum is taken over all forests of the graph~$G$ that consist of~$k$ trees such that each tree $T_j$ in it contains exactly one vertex from the set $\{i_1,\dots i_k\}$.
\end{ter}

\section {Calculation of multidimensional Gaussian sums over a finite field}
Let us explicitly evaluate the expression
$
 \sum_{x_V  \in \mathbb{F}_q^{V} } \chi_1(x_V^t B x_V ).
$
Here $B$ is an arbitrary symmetric matrix of the dimension $|V|\times|V|$ with elements in the field~${\mathbb F}_q$, where $q$ is odd. In this section for convenience we identify the finite set~$V$ with the starting point of some part of the set of natural numbers. As is well known, any symmetric matrix of the rank~$r$ has a nondegenerate principal minor of the order~$r$. This fact is used in the following lemma.

\begin{lem}
\label{lem:4} Let $q=p^d$ with odd prime $p$, $\operatorname{rank}\, B = r$. Then
$$
\sum_{x_V  \in \mathbb{F}_q^V } \chi_1(x_V^t B x_V )  =
q^{|V|} \eta(\det B_r) \left[ {\frac{g(q)}{q}} \right]^r,
$$
where $g(q)$ obeys formula~(\ref{gq}), $\det B_r$ is an arbitrary nonzero principal minor of the order~$r$.
\end{lem}

Before proving Lemma~\ref{lem:4} let us recall some properties of its one-dimensional variant~\cite{ai,lidlNider}.

The value $\sum_{k\in\mathbb F^*_q} \eta ( k )\chi_1(k t) $ is called the Gaussian sum. According to elementary properties of quadratic residues of a finite field, the Gaussian sum vanishes with~${t=0}$. Otherwise, by summing it up with equality~(\ref{eq:5}) and performing certain elementary transformations, we find its value as
\[
\sum_{k\in\mathbb F_q} \chi_1(k^2 t) .
\]

For $t\ne 0$ we can change the variable $x=kt$ in the initial definition of the Gaussian sum and thus find its value as $\eta(t) g(q)$, where
$$
g(q)=\sum_{x\in\mathbb F^*_q} \eta ( x )\chi_1(x).
$$
It is well known that~\cite[theorem~5.15]{lidlNider} the value $g(q)$ obeys formula~(\ref{gq}). As a result we get
\begin{equation}
\label{onedim}
\sum_{k\in\mathbb F_q} \chi_1(k^2 t)=\left\{\begin{array}{ll}
q, & \mbox{if $t=0$,}\\
\eta(t) g(q), & \mbox{if $t\ne 0$.}
\end{array}
\right.
\end{equation}

\textbf{Proof of Lemma~\ref{lem:4}:}\quad 
1. Let us first consider a particular case, namely, let $B$ be diagonal. This means that the left-hand side of the desired equality is representable as follows:
$$
\sum_{x \in \mathbb{F}_q^V}  \chi_1\left(
\sum_{i=1}^{|V|} b_{ii} x_{i}^2 \right).
$$

Let us represent the latter value as the sum of two terms. Namely, let the first term represent the sum of elements with coefficients $b_{ii}$ which are either equal to zero or not (without loss of generality we assume that nonzero coefficients occupy the first places). We get
$$\prod_{i = r + 1}^{|V|} \left( \sum_{x_{i} \in \mathbb F_q }
\chi_1( b_{ii} x_i^2 )   \right) \prod_{i=1}^r \left(\sum_{x_{i} \in \mathbb F_q } \chi_1(b_{ii} x_{i }^2 )
\right) .
$$

Let us apply formula~(\ref{onedim}). Represent the latter expression as
$$
q^{|V|-r} \prod_{i = 1}^r \left( \eta ( b_{ii} ) g(q)
\right) = q^{|V|} \eta \left(
{\det B_r} \right)\left[ \frac{g(q)}{q} \right]^r.
$$

2. The general case can be reduced to the diagonal one. Any nondegenerate symmetric $r\times r$-matrix $B_r$ is reducible to the diagonal form over the field~$\mathbb F_q$. This means that there exists a nondegenerate matrix $Q_r$ of the same dimension such that
$$
  Q_r^t B_r Q_r  = \Lambda,\quad
  \mbox{where}~~\Lambda  = \left( \begin{array}{ccc}
   \lambda _1  &  \cdots  & 0  \\
    \vdots  &  \ddots  &  \vdots   \\
   0 &  \cdots  & \lambda _r    \end{array}  \right)
   ~~\mbox{and $\lambda _i  \ne 0$ for all~$i$.}
$$
Determinants of these matrices satisfy the correlation
\begin{equation}
\label{det}
  \det \Lambda  = (\det Q_r )^2 \det B_r\,.
\end{equation}

In what follows without loss of generality we assume that the matrix $B_r$ is formed by first~$r$ rows and columns of the matrix~$B$. Let us now construct two more nondegenerate matrices; let their dimension equal $|V|\times|V|$. We get
$$
  \widehat Q = \left( \begin{array}{cccc}
   Q_r  & 0 &  \cdots  & 0  \\
   0 & 1 &  \cdots  & 0  \\
    \vdots  &  \vdots  &  \ddots  &  \vdots   \\
   0 & 0 &  \cdots  & 1  \end{array}  \right),\qquad
\widehat Q^t  = \left( {\begin{array}{cccc}
   {Q_r^t } & 0 &  \cdots  & 0  \\
   0 & 1 &  \cdots  & 0  \\
    \vdots  &  \vdots  &  \ddots  &  \vdots   \\
   0 & 0 &  \cdots  & 1   \end{array} } \right).
$$
Let $D = \widehat Q^t B \widehat Q$. Note that properties of the rank imply that $\operatorname{rank}\, D=r$. Here the matrix~$D$ takes the form
$$
  D = \left( \begin{array}{cc}
   \Lambda  & D_{12}   \\
   D_{21}  & D_{22}    \end{array}  \right),
$$
where $D_{12}, D_{21}$ and $D_{22}$ are some (rectangular) matrices obtained as a result of the transformation.

The biunique change of variables $x_V=Q y_V$ gives
$$
\sum_{x_V  \in \mathbb{F}_q^V } \chi_1(x_V^t B x_V )  =
\sum_{y_V  \in \mathbb{F}_q^V } \chi_1(y_V^t \widehat Q^t B
\widehat Q y_V ) = \sum_{y_V  \in \mathbb{F}_q^V } \chi_1(y_V^t D y_V ).
$$

Applying the formula proved in item~1 for the matrix~$D$ of the rank~$r$, we get
$$
\sum_{y_V  \in \mathbb{F}_q^V } \chi_1(y_V^t D y_V )  =
q^{|V|} \eta ( \det D_r ) \left[ \frac{g(q)}{q}
\right]^r;
$$
here $\det D_r$ is a nonzero principal minor of the matrix $D$ of the order $r$ (we use~$\det \Lambda$ as this minor). On the other hand, we know (see~(\ref{det}) that it obeys the formula
$$
\det D_r =
  \left( \det Q_r \right)^2 \det (B_r).
$$
Using properties of the quadratic character $\eta$, we get
$$
  \eta \left( {\det D_r} \right) =
  \eta \left( (\det B_r) \left( \det Q_r \right)^2  \right) =
  \eta \left( \det B_r \right).
$$
This means that
$$
\sum_{x_V  \in \mathbb{F}_q^V } {\chi_1(x_V^t B x_V ) }  = q^{|V|}
\eta(\det B_r) \left[ \frac{g(q)}{q} \right]^r.
$$
$\square$

\section{Proof of the main theorem}

In this section we prove the main theorem treating it as a corollary of results obtained in three previous sections.

\textbf{Proof of Theorem~\ref{ter:1}:} 
Let us transform formula~(\ref{eq:7}) given in the Key lemma~\ref{lem:2}. Applying Lemma~\ref{lem:3} to the inner sum of characters of this expression, we represent this sum as follows:
$$ \sum_{ x_V  \in \mathbb{F}_q^V } \chi_1\left(\sum_{e\in E(G)} (x_{i(e)} - x_{f(e)} )^2 \alpha _e
\right)  = \sum_{  x_V  \in \mathbb{F}_q^V } \chi_1(x_V^t L(G,\alpha) x_V) . $$
By using Lemma~\ref{lem:4} we deduce
\begin{equation}
\label {main}
 F_G (q) = \sum_{\alpha  \in \left( \mathbb{F}_q^*\right)^{E(G)}} \eta \left( \det
L_r(G,\alpha) \right)\left[ \frac{g(q)}{q}
\right]^r,
\end{equation}
where $r$ is the rank of the Laplacian matrix~$L(G,\alpha)$ and $\det L_r(G,\alpha)$ is its nonzero principal minor of the order~$r$. Finally, using Theorem~\ref{ter:3} for $\det L_r(G,\alpha)$, we come to the assertion of Theorem~\ref{ter:1}. 
$\square$

\section{Simplifications and applications of the main theorem}
\subsection{A simplified variant of the main formula and the Tutte conjecture}
In this item we simplify formula~(\ref{eq:3}) given in the main theorem and briefly discuss various related conjectures.

According to Tutte's 3-flow conjecture, for any 4-edge-connected graph~$G$ it holds ${F_G(3)>0}$. This conjecture remained completely inaccessible for a long time, moreover, nothing was known even about $k$- edge-connected graphs. Recently Carsten Thomassen~\cite{thomassen} has succeeded in proving the existence of a 3-flow for 8-edgeconnected graphs. Sometime later this result was strengthened by a group of authors~\cite{combB}, namely, it was proved that the same is true for 6-edge-connected graphs (in this connection see the essay~\cite{lovasz}). However, there are still no arguments in favor of the most intriguing Tutte conjecture, according to which for \textit{any} connected graph~$G$ without bridges it holds $F_G(5)>0$. As far as we know, there are no methods that take into account the specificity of the prime odd number~$q$ or its degree when calculating $F_G(q)$.

Evidently, one can graduate the sum mentioned in the assertion of the main theorem with respect to values~$r^*(G,\alpha)$, namely,
$$
F_G (q) = \sum_{r=0}^{|V|-1} \left[ \frac{g(q)}{q}
\right]^{r} S(r,q),\quad\mbox{ where }
S(r,q)=\sum_{\begin{subarray}{e }\alpha:\,  \alpha \in \left(
\mathbb{F}_q^*\right)^{E(G)},\\\phantom{\alpha:}\, r^*(G,\alpha)=r\end{subarray}} {\eta
\left(s(\alpha,G/W^*)\right)}.
$$

\begin{utv}
\label{utv:1}
For any odd~$r$ it holds $S(r,q)=0$.
\end{utv}

\textbf{Proof:}
Really, let $\gamma$ be an arbitrary element of the field ${\mathbb F}_q$ such that $\eta(\gamma)=-1$. Let us associate an arbitrary collection $\alpha_e, e\in E(G)$ with a collection $\beta_e, e\in E(G)$, where
$\alpha_e=\gamma \beta_e$. Evidently, this correspondence is biunique. The number of vertices in the contracted graph $G/W^*$ equals ${|V|-|W|+1}$. Therefore in any tree $T$ in this graph the number of edges equals $|V|-|W|$, consequently,
$\prod_{e\in E(T)} \beta_e =\prod_{e\in E(T)} \alpha_e \gamma^{|V|-|W|}$. Thus, we get $s(\alpha,G/W)=s(\beta,G/W) \gamma^{|V|-|W|}$. In view of the multiplicativity of the symbol $\eta$ we obtain
$$
\sum_{\begin{subarray}{e }\alpha:\,  \alpha \in \left(
\mathbb{F}_q^*\right)^{E(G)},\\\phantom{\alpha:}\, r^*(G,\alpha)=r\end{subarray}} {\eta
\left(s(\alpha,G/W^*)\right)} =\eta(\gamma)^r
\sum_{\begin{subarray}{e }\beta:\, \beta \in \left(
\mathbb{F}_q^*\right)^{E(G)},\\\phantom{\beta:}\, r^*(G,\beta)=r\end{subarray}} {\eta
\left(s(\beta,G/W^*)\right)} =-S(r,q),
$$
consequently, in this case $S(r,q)=-S(r,q)=0$.
$\square$

Therefore, the main formula obtained in this paper is representable as follows:
\begin{equation}
\label{eq:4MT} F_G (q) = \sum_{\begin{subarray}{e} r:\, r \bmod 2 =
0, \\ \phantom{r:}\,r\leqslant |V|-1,
\end{subarray}} \left[ \frac{g(q)}{q} \right]^{r} S(r,q).
\end{equation}

In the remaining part of the paper we discuss various conjectures connected with
values $S(r,q)$ when $q=5$. Note that coefficients of the linear combination $S(r,q)$ in formula~(\ref{eq:4MT}) with $q=5$ are always positive ($q=5$ is the first prime number with this property). Really, let $r=2 r'$, then with prime $q$ such that $q \bmod 4 = 1$ we get
$\left[ \frac{g(q)}{q} \right]^{2r'}=\frac{1}{q^{r'}}$.

Consequently, the flow polynomial $F_G(q)$ is automatically positive with $q=5$, if at least one value of $r$ in formula~(\ref{eq:4MT}) ensures that $S(r,5)> 0$, while with the rest values of $r$ it holds $S(r,5)\geqslant 0$. Unfortunately, (as we see in examples given in subsection~7.3) this is not necessarily true for an arbitrary graph. In the next item we give some facts that are useful in calculations.

\subsection{The rank of the Laplacian matrix over a finite field}
Here we study some easy properties of the Laplacian matrix which enable us to estimate the value $r^*(G,\alpha)$. Let us first pay attention to one simple property of formula~(\ref{main}).
\begin{utv}
\label{utv:2}
If a graph~$G$ contains at least one nonmultiple edge (in particular, if the graph $G$ is simple and nontrivial), then formula~(\ref{eq:4MT}) contains no terms that correspond to $r=0$.
\end{utv}
\textbf{Proof:}
Really, the rank of the matrix differs from zero if the matrix has at least one nonzero element. If the graph $G$ has at least one nonmultiple edge $e = (i(e), f(e))$, then the element of the matrix $L(G,\alpha)$ with indices $i(e)$, $f(e)$ equals $\alpha_e \ne 0$.
$\square$

In what follows we assume that $q\equiv p$ is an odd prime number. The easy property given below allows us to estimate possible values of $|W^*|$ (or, equivalently, to estimate $r^*(G,\alpha)$), using the known integer value of the Laplacian determinant.
\begin{lem}
\label{linal}
Let $p$ be a prime number and let $B$ be an $n\times n$ integer matrix of the rank~$r$ considered as a matrix over the residue field modulo~$p$. Then the determinant of any minor of the matrix of the order $r+i$, $i\in\mathbb N$, $i\leqslant n-r$, considered as an integer number, is a multiple of~$p^i$.
\end{lem}
\textbf{Proof:}
If the rank of the matrix~$B$ over the field~$\mathbb F_p$ equals~$r$, then any its minor~$M$ of the order $r+i$ contains no more than $r$ independent basic rows over the field~$\mathbb F_p$. Then each of all the rest rows in this minor is representable as a linear combination of these rows. The addition to some row of a linear combination of the rest ones does not affect the value of the determinant over any field. Consequently, when evaluating the determinant we can assume that all elements in these row equal zero modulo~$p$, i.e., they are multiples of $p$. Then the determinant of the minor~$M$ (over the field~$\mathbb R$) is a multiple of $p^j$, where $j$ is the number of the mentioned rows.
$\square$

\begin{sled}
Let $G$ be a connected multigraph without loops and let $q=p$ be an odd prime number. Assume that the determinant of the matrix $L(G,\alpha)$ is not a multiple of $p^{k+1}$, $k\in\mathbb N$. Then $|W^*|\leqslant k$ and $r^*(G,\alpha)\geqslant |V|-k$, respectively.
\end{sled}
\textbf{Proof:}
Really, by definition the principal minor of the matrix $L(G,\alpha)$ that corresponds to vertices in $V-W^*$ differs from zero and $W^*$ is the minimal cardinality set with such a property. Therefore $\operatorname{rank} L(G,\alpha)=|V|-|W^*|$. Consequently, by Lemma~\ref{linal}, if $|W^*|>k$, then $L(G,\alpha)$ would be a multiple of~$p^{k+1}$.
$\square$

\subsection{Examples of calculations for $q=5$}
We have considered various graphs, which are basic in Tutte's conjecture. Evidently, if the degree of some vertex $v$ of the graph $G$ equals~2, then the evaluation of the flow polynomial $F_G(q)$ in this graph is reduced to the evaluation of $F_{G'}(q)$, where $G'$ is obtained from $G$ by contracting the vertex $v$ and an adjacent one.
For this reason we have considered only those graphs where the degree of vertices is not less than three.

It is well known (see \cite{lovasz}) that with $q\geqslant 3$ an example of a graph without bridges such that $F_G(q) < 0$ and the sum $|E|+|V|$ takes on the minimal value is a simple cubic graph (provided that such an example exists). Thus, for $q=3$ the graph $K_4$ serves as such an example, and for $q=4$ the Petersen graph does. According to the Tutte conjecture, for $q=5$ such an example does not exist. Therefore, for proving this conjecture it suffices to consider simple cubic graphs.

We have considered $S(r,5)$ for all such graphs with no more than 10 vertices. As appeared, all such graphs with less than 10 vertices satisfy the condition $S(r,5)\geqslant 0$ for all~$r$. However, we have found several graphs with 10 vertices for which $S(6,5)< 0$. In particular, for the Petersen graph we have $S(2,5)=S(4,5)=0$, $S(6,5) =  -384$, $S(8,5) = 151\,920$, and $F_G(5) = S(6,5)/125+S(8,5)/625=240$.

Moreover, one can find such examples even for simple noncubic graphs with lesser numbers of vertices. Thus, for the graph $K_{5}$ with one added vertex of the degree~3 we have obtained the following values: $S(2,5) =  -180$, $S(4,5) = 513\,300$, and $F_G(5) =S(2,5)/5+S(4,5)/25=20\,496$. One may think that the last nonzero term of the sequence $S(i,5), i=2,4,\ldots$, is positive for any simple graph of all degrees with more than two vertices. However, this is not true. For the graph $K_{3,4}$ it holds $S(2,5) = 612$, $S(4,5)=244\,860$, $S(6,5)=-8\,100$, and $F_G(5)=S(2,5)/5+S(4,5)/25+S(6,5)/125=9\,852$.

Therefore, the reason, for which for any connected graph without bridges (in accordance with the Tutte conjecture) the sum of Legendre symbols with coefficients $\frac15,\frac{1}{25}$, etc. given in Theorem~\ref{ter:1} equals a positive (integer) number, remains unknown.

\end{document}